\documentclass [a4paper] {article}

\usepackage[utf8]{inputenc}

\usepackage[french]{babel}

\usepackage{amssymb,amsmath, amsfonts}

\usepackage [dvips] {graphicx,color}

\usepackage{amsthm}

\usepackage{MnSymbol}

\newtheorem{prop} {Proposition} 
\newtheorem{lm} [prop]{Lemme} 
\newtheorem{thm} [prop] {Théorème} 
\newtheorem{cor} [prop] {Corollaire}

\theoremstyle{definition}
\newtheorem{df}{Définition} 
\newtheorem*{df*}{Définition}

\theoremstyle{remark}
\newtheorem{rmq}{Remarque} 
\newtheorem{example}{Exemple} 
\newtheorem{exm}[example]{Exemple}

\author{Stéphane \textsc{Dugowson}
\footnote{Laboratoire LISMMA-QUARTZ, Institut Supérieur de Mécanique de Paris. Email :  stephane.dugowson@supmeca.fr}
}

\title {Structure connective des relations multiples}

\begin{document}

\maketitle

\noindent\textbf{Résumé.} L'objet de cet article est d'abord de définir la structure connective sur un ensemble $I$ de toute relation multiple portant sur une famille d'ensembles indexée par $I$, une telle relation étant vue comme exprimant une compatibilité entre les états de différents systèmes, de sorte qu'une compatibilité totale exprime en fait une absence d'interaction. Nous démontrons alors un \og théorème de Brunn\fg\, pour les relations multiples, à savoir le fait que toute structure connective est celle d'une telle relation.

\noindent\textbf{Mots clés.} Connectivité. Espaces connectifs. Relations. Algèbre relationnelle. Borroméen. Brunn. 

\mbox{}

\noindent\textbf{Abstract.} The prime purpose of this paper is to define the connectivity structure, on a set $I$, of any multiple relation defined on a family of sets indexed by $I$, such a relation expressing compatibility between the states of different systems (thus a full compatibility indicates absence of any connection). We then demonstrate a "Brunn's theorem" for those multiple relations, that is the fact that every connectivity structure is the connectivity structure of such a relation.


\mbox{}


%
%

\mbox{}

Les relations multiples considérées ici sont d'arité quelconque variable, gé\-né\-ra\-lement infinie, puisque les ensembles entre les éléments desquels elles portent sont indexés par un sous-ensemble d'un ensemble $I$ quelconque\footnote{Malgré cela, on ne fera pas a priori l'hypothèse de l'axiome du choix. Du reste, pour tous les exemples auxquels nous avons songé en pratique, les ensembles de la famille considérée ont des éléments repérables qui peuvent être choisis directement, de sorte que le produit de ces ensembles est assurément non vide.} fixé. Intuitivement, les ensembles ainsi indexés représentent autant de systèmes en interaction mutuelle, leurs éléments représentent les états de ces systèmes, et une relation multiple exprime la compatibilité de certains de ces états, donc certaines contraintes mutuelles entre les systèmes considérés. Afin de définir la structure connective sur $I$ d'une telle relation, nous commençons par préciser dans la première partie, consacrée à un certain monoïde commutatif et idempotent constitué de ces relations multiples, certaines notations, définitions et résultats relatifs à de telles relations. Cette partie n'a aucune prétention à l'originalité, on en trouve sans doute la substance dans les cours d'\emph{algèbre relationnelle} --- même si l'on s'y limite généralement aux arités finies --- et pour la plupart les résultats donnés sont intuitivement évidents et presque toujours très faciles à démontrer. La deuxième partie, plus originale, où le résultat le plus intéressant est sans doute que l'union de parties dites détachables n'est pas nécessairement elle-même détachable, prépare à la troisième, où se trouve définie la structure connective d'une relation multiple. On démontre finalement un \og théorème de Brunn\fg pour les relations multiples, à savoir le fait que toute structure connective est celle d'une relation multiple.

\mbox{}

Dans tout l'article, $I$ désigne un ensemble, et $\mathcal{E}=(E_i)_{i\in I}$ une famille d'ensembles non vides indexée par $I$. On notera $\vert \mathcal{E}\vert$ l'ensemble $\cup_{i\in I} E_i$. Pour tout ensemble $A$, on désigne par $\mathcal{P}A$ l'ensemble des parties de $A$.

\section{Le monoïde commutatif $(\mathcal{R}_\mathcal{E}, \bowtie, 1)$}

\subsection{Graphes triviaux et familles}

\begin{df} Pour toute partie $J$ de $I$, on appelle \emph{graphe total} ou \emph{graphe trivial} sur $J$ et on note $Z_J$ le produit cartésien
\[Z_J=\prod_{j\in J} E_j.\] 
\end{df}

Les élément de $Z_J$ seront appelés des \emph{$J$-familles}. Cette dénomination est sans ambiguïté car les $Z_J$ sont deux à deux disjoints. L'ensemble de toutes les \emph{familles dans $\mathcal{E}$} sera noté $\mathcal{Z}$ :
\[\mathcal{Z}=\bigcup_{J\in\mathcal{P}(I)} Z_J.\]

Pour $x\in\mathcal{Z}$, nous noterons $J_x$ le \emph{domaine} de $x$, c'est-à-dire l'unique partie $J_x$ de $I$ telle que $x\in Z_{J_x}$.

\begin{rmq} Toute $J$-famille $x\in Z_J$ s'identifie à une application $x:J\to \vert \mathcal{E}\vert$ vérifiant la condition suivante :
\[\forall j\in J, x_j\in E_j,\] où $x_j=x(j)$ désigne l'image de $j$ par cette application.
\end{rmq}

\begin{rmq} [Cas où $J=\emptyset$] Le graphe total sur $\emptyset\subset I$ est réduit à un singleton, dont par convention l'unique élément sera noté $\bullet$ :
\[Z_\emptyset=\{\bullet\}.\] Il y donc une unique $\emptyset$-famille, à savoir $\bullet$.
\end{rmq}

\subsection{Relations multiples}

\begin{df} Une \emph{relation multiple} $R$ dans $\mathcal{E}$ est un couple $(J,G)$ constitué
\begin{itemize}
\item d'une partie $J\subset I$ appelée \emph{domaine} de $R$,
\item d'une partie $G\subset Z_J$ appelée le \emph{graphe} de $R$.
\end{itemize}
\end{df}

Dans la suite, les relations multiples dans $\mathcal{E}$ seront également appelées plus simplement des \emph{$\mathcal{E}$-relations}, ou plus simplement encore des \emph{relations}. 

Une relation $R$ étant donnée, nous désignerons parfois par $J_R$ son domaine, et par $G_R$ son graphe. Les relations de domaine $J$ pourront être appelées des $J$-relations. Dans le cas où $card(J)=2$, on parlera de relations binaires.

Les $J_R$-familles appartenant à $G_R$ seront dites \emph{compatibles pour la relation $R$}, ou encore \emph{$R$-compatibles}. Par abus d'écriture, on écrira souvent $x\in R$ au lieu de $x\in G_R$ pour exprimer qu'une famille $x\in Z_{J_R}$ est $R$-compatible.

L'ensemble des relations multiples dans $\mathcal{E}$ sera noté $\mathcal{R}_\mathcal{E}$ ou plus simplement, puiqu'ici nous considérons que $\mathcal{E}$ est fixé, $\mathcal{R}$.
L'ensemble des relations multiples de domaine $J\subset I$ sera noté $\mathcal{R}_J$, de sorte que
\[\mathcal{R}=\bigcup_{J\subset I} \mathcal{R}_J.\]

Pour tout $J\subset I$, les $J$-relations sont ordonnées par l'inclusion de leurs graphes. Étant données deux $J$-relations $R$ et $S$, nous écrirons $R\subset_J S$, ou simplement $R\subset S$, pour exprimer le fait que $G_R\subset G_S$. 

\subsubsection{Relations nulles et relations triviales}

La $J$-relation minimale, notée $0_J$, est celle de graphe vide :
\[0_J=(J,\emptyset),\]
tandis que la $J$-relation maximale, notée $1_J$, est celle de graphe total
\[1_J=(J,Z_J).\]

Les relations de la forme $0_J$ seront dites \emph{nulles}, celles de la forme $1_J$ seront dites \emph{triviales} (ou \emph{totales}).

Si l'ensemble $J$ est fini, ou si l'on admet l'axiome du choix, alors $Z_J\neq\emptyset$, de sorte que $0_J\neq 1_J$. Ceci est vrai, bien que ce ne soit pas très intuitif, en particulier pour $J=\emptyset$, auquel cas le graphe de $0_\emptyset$ est vide tandis que celui de $1_\emptyset$ est $Z_\emptyset=\{\bullet\}$.

Dans la suite, on notera également $\mathbf{1}$ cette dernière relation \og sur aucun ensemble\fg  :
\[\mathbf{1}=1_\emptyset=(\emptyset, \{\bullet\}).\]

On posera également

\[\overline{\mathbf{1}}=1_I,\]
et
\[\mathbf{0}=0_I.\]

\subsection{Restrictions}

Dans cette section, on considère deux parties $J$ et $K$ de $I$ telles que $J\subset K\subset I$.

\subsubsection{Restrictions de familles}

\begin{df}  On appelle  \emph{restriction (de K) à J}, l'application $\rho_{(K,J)}:Z_K\to Z_J$ définie pour tout $x\in Z_K$ par
\[\rho_{(K,J)}(x)=x\circ (J\hookrightarrow  K),\] où $x$ est vu comme application $x:K\to\vert \mathcal{E}\vert$ et
où $(J\hookrightarrow  K)$ désigne l'injection canonique de $J$ dans $K$.
Autrement dit, 
\[\rho_{(K,J)} ((x_k)_{k\in K})=(x_k)_{k\in J}.\] 
\end{df}

L'image $(x_k)_{k\in J}$ d'un élément $(x_k)_{k\in K}$ de $Z_K$ par $\rho_{(K,J)}$ est la \emph{restriction à $J$} de $(x_k)_{k\in K}$. Bien entendu, pour $K=J$, on a 
\[\rho_{(J,J)}=id_{Z_J}.\]

\begin{exm}  La restriction à $J=\emptyset$ d'une $K$-famille quelconque $x$ est $\bullet$.
\end{exm}

\begin{df} Soit $J$, $K$ et $L$ trois parties de $I$ telles que
\[L\subset J\cap K.\] Soit $x\in Z_J$ et $y\in Z_K$. On dit que \emph{$x$ et $y$ coïncident sur $L$} si
\[\rho_{(J,L)}(x)=\rho_{(K,L)}(y).\]
\end{df}

\subsubsection{Restriction des relations}

Supposant toujours $J\subset K\subset I$, la notion de restriction, définie sur les familles $(x_k)_{k\in K}$, s'étend évidemment aux relations elles-mêmes :

\begin{df} L'application $\mathcal{R}_K\to \mathcal{R}_J$ qui à toute relation multiple $R=(K,G)$ avec $G\subset Z_K$ associe la relation $(J,H)$ avec 
\[H=\rho_{(K,J)}(G)=\{\rho_{(K,J)}(x), x\in G\}\subset Z_J\]
sera encore notée $\rho_{(K,J)}$ et encore appelée \emph{restriction (de $K$) à $J$}. On a ainsi
\[\rho_{(K,J)}(K,G)=(J, \rho_{(K,J)}(G))\in\mathcal{R}_J.\]
\end{df}

\begin{exm} Dans le cas où $J=\emptyset$, la restriction à $J$ de la $K$-relation vide $0_K)$ est $0_\emptyset$, tandis que la restriction à $J$  d'une relation non vide quelconque est $\mathbf{1}=1_\emptyset$. 
\end{exm}

\subsubsection{Composition des restrictions}

Qu'elles portent sur les familles ou sur les relations, les restrictions se composent évidemment selon 
\[\rho_{(K,L)}\circ \rho_{(J,K)}=\rho_{(J,L)},\] où l'on a supposé $I\supset J\supset K\supset L$.

\subsubsection{Notation}
Plutôt que de noter $(\rho_{(K,L)}\circ \rho_{(J,K)})(x)$ la restriction successive d'une $J$-famille $x$ à $K$ puis à $L$, il serait plus commode de l'écrire
\[x\rho_{(J,K)}\rho_{(K,L)}.\] Nous n'adopterons pas ici cette remise à l'endroit des notations de composition, mais les notations usuelles pour la restriction, qui n'explicitent pas le domaine de départ mais uniquement le domaine d'arrivée de la restriction, seront très utiles, et nous poserons ainsi pour toute $J$-famille $x$ et toute partie $K\subset J$ :
\[\rho_{(J,K)}(x)=x_{\vert K}.\] De même pour les restrictions de relations :
\[\rho_{(J,K)}(R)=R_{\vert K}.\] Avec ces notations, la composition des restrictions s'écrit simplement
\[{x_{\vert K}}_{\vert L}=x_{L}.\]

\subsection{Sommes de familles}

\subsubsection{Familles compatibles entre elles}

\begin{df} Pour tout couple de familles $(x,y)\in \mathcal{Z}^2$, on dit que \emph{$x$ et $y$ sont compatibles entre elles}, et on note $x\diamond y$, si $x$ et $y$ coïncident sur l'intersection de leurs domaines, autrement dit si
\[x_{\vert J_x\cap J_y}=y_{\vert J_x\cap J_y}.\]
\end{df}

La relation $\diamond$ est évidemment réflexive et symétrique sur $\mathcal{Z}$.

\subsubsection{Somme de deux familles compatibles}

\begin{df} On définit ainsi une opération binaire partielle sur $\mathcal{Z}$ : pour tout couple de familles $(x,y)\in \mathcal{Z}^2$ tel que $x\diamond y$, on note $x+y$ la famille de domaine 
\[J_{x+y}=J_x\cup J_y,\]
et telle que 
\[\forall j\in J_x, (x+y)_j=x_j \quad\mathrm{et}\quad \forall j\in J_y, (x+y)_j=y_j.\]
\end{df}

Autrement dit, $x+y$ est caractérisée par son domaine $J_x\cup J_y$ et par le fait que 
\[(x+y)_{\vert J_x}=x\quad\mathrm{et}\quad(x+y)_{\vert J_y}=y.\]

On vérifie immédiatement la proposition suivante.

\begin{prop}[Elément neutre, idempotence et commutativité] L'opération $+$ vérifie les trois propriétés suivantes :
\begin{itemize}
\item pour tout $x\in\mathcal{Z}$, $x+\bullet=x$,
\item pour tout $x\in\mathcal{Z}$, $x+x=x$,
\item si $x$ et $y$ sont deux familles compatibles, alors $x+y=y+x$.
\end{itemize}
\end{prop}

\subsubsection{Somme de restrictions}

\begin{prop}\label{prop somme de deux restrictions} Pour tout $x\in\mathcal{Z}$ et tout couple $(K,L)$ de parties de $J_x$ on a 
\[x_{\vert K}+ x_{\vert L}=x_{\vert K\cup L}.\]  En particulier, si $K\cup L=J_x$, on a $x=x_{\vert K}+ x_{\vert L}$.
\end{prop}
\paragraph{Preuve.} On a $x_K\diamond x_L$ puisque ${x_{\vert K}}_{\vert K\cap L}={x}_{\vert K\cap L}={x_{\vert L}}_{\vert K\cap L}$. Et $x_{\vert K}+ x_{\vert L}$ coïncide trivialement avec $x$ sur $K\cup L$, d'où l'égalité annoncée. 
\begin{flushright}$\square$\end{flushright} 

\subsubsection{Restriction de somme}

\begin{prop} Soient $x$ et $y$ deux familles compatibles dans $\mathcal{E}$, et soit $L\subset J_x\cup J_y$. On a
\[(x+y)_{\vert L}=x_{\vert J_x\cap L}+y_{\vert J_y\cap L}.\]
\end{prop}
\paragraph{Preuve.} On applique la proposition précédente au recouvrement de $L$ par $J_x\cap L$ et $J_y\cap L$, d'où
\[(x+y)_{\vert L}=((x+y)_{\vert L})_{\vert J_x\cap L}+((x+y)_{\vert L})_{\vert J_y\cap L}.\]
Mais $((x+y)_{\vert L})_{\vert J_x\cap L}=(x+y)_{\vert J_x\cap L}=((x+y)_{\vert J_x})_{\vert J_x\cap L}=x_{\vert J_x\cap L}$.
De même, $((x+y)_{\vert L})_{\vert J_y\cap L}=y_{\vert J_y\cap L}$. D'où le résultat.
\begin{flushright}$\square$\end{flushright}

\subsubsection{Familles finies de familles}

\begin{prop} Soient $x$, $y$ et $z$ trois familles dans $\mathcal{E}$. Si ces familles sont deux à deux compatibles :
$x\diamond y$, $y\diamond z$, et $z\diamond x$, alors
$(x+y)\diamond z$.
\end{prop}
\paragraph{Preuve.} En effet, $(x+y)_{\vert (J_x\cup J_y)\cap J_z}=x_{\vert J_x\cap J_z}+y_{\vert J_y\cap J_z}=z_{\vert J_x\cap J_z}+z_{\vert J_y\cap J_z}=z_{\vert (J_x\cup J_y)\cap J_z}$.
\begin{flushright}$\square$\end{flushright}

\begin{prop}[Associativité] L'opération binaire partielle $+$ est associative sur $\mathcal{Z}$ au sens où pour tout triplet $(x,y,z)\in\mathcal{Z}^3$,
\[x\diamond y\diamond z\diamond x \Rightarrow (x+y)+z=x+(y+z).\]
\end{prop}
\paragraph{Preuve.} L'existence des sommes considérées est assurée par la proposition précédente, et leur égalité se vérifie immédiatement.
\begin{flushright}$\square$\end{flushright} 

Du fait de l'associativité et de la commutativité de l'opération $+$, nous pourrons parler de la somme $\sum_{\lambda\in\Lambda}x_\lambda$ d'une famille finie\footnote{Et la notion s'étend sans difficulté à une famille quelconque de familles deux à deux compatibles.} quelconque $(x_\lambda)_{\lambda\in\Lambda}$ de familles deux à deux compatibles dans $\mathcal{E}$. En particulier, la proposition \ref{prop somme de deux restrictions} se généralise facilement :

\begin{prop} Pour tout $x\in\mathcal{Z}$ et tout recouvrement fini $(J_\lambda)_{\lambda\in\Lambda}$ de $J_x$, on a 
\[x=\sum_{\lambda\in\Lambda}x_{\vert J_\lambda}.\]
\end{prop}

\subsection{Produit de relations}

Pour désigner le produit des relations considéré ici, nous reprenons la notation $\bowtie$, usuelle en algèbre relationnelle pour désigner  la \emph{jointure} de deux relations.

\subsubsection{Définition du monoïde commutatif $(\mathcal{R}_\mathcal{E}, \bowtie, 1)$}

\begin{df}[Produit de deux relations multiples] Étant données $R=(J_R,G_R)$ et $S=(J_S,G_S)$ deux relations multiples dans $\mathcal{E}$, on définit leur produit $T=R\bowtie S$ de la façon suivante :
\begin{itemize}
\item le domaine $J_T$ de $T$ est l'union $J_T=J_R\cup J_S$,
\item le graphe $G_T\subset Z_{J_T}$ de $T$ est défini par 
\[G_T=\{x\in Z_{J_T}, \rho_{(J_T,J_R)}(x)\in G_R\,\,\mathrm{et}\,\,\rho_{(J_T,J_S)}(x)\in G_S\}.\] Autrement dit, 
\[(J_R,G_R)\bowtie (J_S,G_S)=(J_R\cup J_S, \rho_{(J_T,J_R)}^{-1}(G_R)\cap \rho_{(J_T,J_S)}^{-1}(G_S)).  \]
\end{itemize}
\end{df}

\begin{prop}\label{prop expression du produit de relations par sommes} Le graphe du produit $R\bowtie S$ de deux relations $R$ et $S$ dans $\mathcal{E}$ est donné par
\[G_{R\bowtie S}=\{r+s, r\in R, s\in S, r\diamond s\}.\]
\end{prop}
\paragraph{Preuve.} Pour $r$ et $s$ comme ci-dessus, on a $\rho_{(J_R\cup J_S,J_R)}(r+s)=(r+s)_{\vert J_R}=r\in G_R$, et de même $\rho_{(J_R\cup J_S,J_S)}(r+s)=(r+s)_{\vert J_S}=s\in G_S$. Réciproquement, pour toute $(J_R\cup J_S)$-famille $x$ telle que $r=\rho_{(J_T,J_R)}(x)\in G_R$ et $s=\rho_{(J_T,J_S)}(x)\in G_S$, on a clairement $r\diamond s$ et $x=r+s$.
\begin{flushright}$\square$\end{flushright} 

\begin{prop} L'opération $\bowtie$ ainsi définie sur l'ensemble $\mathcal{R}$ des relations multiples dans $\mathcal{E}$ est associative, commutative, idempotente et admet pour élément neutre la \og relation pleine sur aucun ensemble\fg\, $\mathbf{1}=1_\emptyset$.
\end{prop}

\paragraph{Preuve.} L'associativité et la commutativité de $\bowtie$ découle  de celles de l'union, de l'intersection et de la composition des opérations de restriction. Pour toute relation $R$, l'égalité
\[R\bowtie R=R\] est immédiate.
Enfin, on a
\[(J,G)
\bowtie \mathbf{1}
=
(J,G)\bowtie (\emptyset,\{\bullet\})
=
(J\cup\emptyset,\rho_{(J,J)}^{-1}(G)\cap \rho_{(J,\emptyset}^{-1}(\{\bullet\}))
=
(J,G\cap Z_J)
=
(J,G).\]
\begin{flushright}$\square$\end{flushright} 

De la définition du produit $R\bowtie S$, on déduit immédiatement que la restriction au domaine de $R$ d'une famille $R\bowtie S$-compatible est nécessairement $R$-compatible : 

\begin{prop}\label{prop restriction de produit inclus} Pour toutes relations $R$ et $S$ dans $\mathcal{E}$, on a 
\[(R\bowtie S)_{\vert J_R}\subset R.\]
\end{prop}

\subsubsection{Exemples}

\paragraph{Produit par $\overline{\mathbf{1}}=1_I$ : prolongement à $I$.}

\begin{df}
Pour tout $J\subset I$ et tout $R\in \mathcal{R}_J$, on appelle \emph{prolongement à $I$} et l'on note $\overline{R}$ la $I$-relation définie par
\[\overline{R}=R\bowtie \overline{\mathbf{1}}. \]
\end{df}

\begin{rmq} La notation $\overline{R}$ est compatible avec celles de $\mathbf{1}$ et de $\overline{\mathbf{1}}$, puisque 
\[\overline{\mathbf{1}}=\mathbf{1}\bowtie\overline{\mathbf{1}}. \]
\end{rmq}

\begin{prop} Pour tout $J\subset I$ et tout $R\in \mathcal{R}_J$, on a
\[\overline{R}=(I,\rho_{(I,J)}^{-1}(G_R))= R\bowtie 1_{\neg J}, \] où $\neg J=I\setminus J$.
\end{prop}

\paragraph{Preuve.} Le domaine de $\overline{R}$ est $I=J\cup I=J\cup \neg J$, et son graphe est
\[G_{\overline{R}}=\rho_{(I,J)}^{-1}(G_R)\cap\rho_{(I,I)}^{-1}(Z_I)= \rho_{(I,J)}^{-1}(G_R).\] D'un autre coté, $\rho_{(I,\neg J)}^{-1}(Z_{\neg J})=Z_I$, de sorte que $\rho_{(I,J)}^{-1}(G_R)=\rho_{(I,J)}^{-1}(G_R)\cap \rho_{(I,\neg J)}^{-1}(Z_{\neg J})$ d'où finalement
\[G_{\overline{R}}=G_{R\bowtie 1_{\neg J}}.\]
\begin{flushright}$\square$\end{flushright} 

\paragraph{Produit par $\mathbf{0}=0_I$.}

Pour toute relation $R\in\mathcal{R}$, il est immédiat que
\[R\bowtie \mathbf{0} = \mathbf{0}.\]

\paragraph{Produit par $0_\emptyset$.}

Pour toute relation $R\in\mathcal{R}$, il est immédiat que $R\bowtie 0_\emptyset$ est la relation vide de même domaine que $R$ :
\[R\bowtie 0_\emptyset = (J_R,\emptyset)=0_{J_R}.\]

\paragraph{Composition de deux relations binaires.}

Supposons que $I$ contienne une partie $J$ ayant trois éléments distincts notés $1$, $2$ et $3$,
\[J=\{1,2,3\}\subset I,\]
 et soient $f:E_1\to E_2$ et $g:E_2\to E_3$ deux relations binaires (par exemple deux applications), de domaines respectifs $J_f=\{1,2\}$ et $J_g=\{2,3\}$. Leur produit est alors défini par
\[f\bowtie g=g\bowtie f =(\{1,2,3\},\{(x,y,z)\in E_1\bowtie E_2\bowtie E_3, y=f(x)\,\,\mathrm{et}\,\, z=g(y)\}),\]
de sorte que 
\[g\circ f=\rho_{(\{1,2,3\},\{1,3\})} (f\bowtie g).\]

\begin{rmq} La non-commutativité de la composition des applications (ou des relations binaires) n'est bien entendu pas contredite par la commutativité du produit des relations, pour lequel l'information des $E_i$ en jeu est contenue dans la donnée du domaine.
\end{rmq}

\subsection{Relations incompatibles}

\begin{df} Deux relations dans $\mathcal{E}$ sont dites incompatibles si leur produit est nul.
\end{df}

La proposition suivante est immédiate.

\begin{prop}Deux relations $R$ et $S$ sont incompatibles si et seulement si pour tout $x\in R$ et tout $y\in S$, on a $x$ et $y$ incompatibles.
\end{prop}

\begin{exm} Toute relation nulle est incompatible avec toute autre relation.
\end{exm}

\subsection{Restriction d'un produit}


%

\begin{prop}\label{prop inclusion restriction produit}
Soient $R$ et $S$ deux relations dans $\mathcal{E}$, et $L$ une partie de $J_R\cup J_S$. On a
\[(R\bowtie S)_{\vert L}\subset R_{\vert L\cap J_R}\bowtie S_{\vert L\cap J_S}.\]
\end{prop}

\paragraph{Preuve.} 

Les deux relations sont comparables, puisqu'elles sont de même domaine $L$. 
Soit maintenant $x=(r+s)_{\vert L}\in (R\bowtie S)_{\vert L}$, avec $r\in R$, $s\in S$ et $r\diamond s$. On a
\[x=x_{\vert L\cap J_R}+x_{\vert L\cap J_S}.\]
Mais 
\[x_{\vert L\cap J_R}=((r+s)_{\vert L})_{\vert L\cap J_R}
= (r+s)_{\vert L\cap J_R}=((r+s)_{\vert J_R})_{\vert L\cap J_R}=r_{\vert L\cap J_R},\] d'où $x_{\vert L\cap J_R}\in R_{\vert L\cap J_R}$. De même,
$x_{\vert L\cap J_S}\in S_{\vert L\cap J_S}$.
On en déduit que $x\in R_{\vert L\cap J_R}\bowtie S_{\vert L\cap J_S}$.
\begin{flushright}$\square$\end{flushright} 

\begin{rmq} La réciproque est fausse en général. Par exemple, pour $I=\{1,2,3,4\}$, $E_i=\mathbf{N}$ pour tout $i\in I$, $J_R=\{1,2,3\}$, $J_S=\{1,2,4\}$, $L=\{1,3,4\}$, $R$ défini par 
\[(r_1,r_2,r_3)\in R \Leftrightarrow r_2=0\]
et $S$ défini par 
\[(s_1,s_2,s_4)\in S \Leftrightarrow s_2\neq 0,\] on a $R$ et $S$ incompatibles de sorte que
\[(R\bowtie S)_{\vert L}=0_L,\] mais par ailleurs
\[R_{\vert L\cap J_R}\bowtie S_{\vert L\cap J_S}=1_{L\cap J_R}\bowtie 1_{L\cap J_S}=1_L.\]
\end{rmq}
 
%

\begin{prop}\label{prop egalite restriction produit} Soient $R$ et $S$ deux relations dans $\mathcal{E}$, et $L$ une partie de $I$ vérifiant
\[J_R\cap J_S\subset L \subset J_R\cup J_S.\]
Alors
\[(R\bowtie S)_{\vert L}= R_{\vert L\cap J_R}\bowtie S_{\vert L\cap J_S}.\] En particulier, si $J_R\cap J_S=\emptyset$, l'égalité ci-dessus est satisfaite pour tout $L \subset J_R\cup J_S$.
\end{prop}

\paragraph{Preuve.} D'après la proposition \ref{prop inclusion restriction produit}, il suffit de prouver l'inclusion du second membre dans le premier. Soit donc $r'+s'\in R_{\vert J_R\cap L}\bowtie S_{\vert J_S\cap L}$, avec $r'\in R_{\vert J_R\cap L}$ et $s'\in S_{\vert J_S\cap L}$ qui coïncident sur $J_R\cap J_S$. 
Puisque $r'\in R_{\vert J_R\cap L}$, il existe $r\in R$ tel que $r_{\vert J_R\cap L}=r'$, et de même il existe $s\in S$ tel que $s_{\vert J_S\cap L}=s'$. 
L'inclusion $J_R\cap J_S\subset L$ entraîne alors d'une part $J_R\cap J_S\subset J_R\cap L$ d'où $r_{\vert J_R\cap J_S}$ $=(r_{\vert J_R\cap L})_{\vert J_R\cap J_S}$  $=r'_{\vert J_R\cap J_S}$, et d'autre part $J_R\cap J_S\subset J_S\cap L$ d'où  $s_{\vert J_R\cap J_S}$ $=s'_{\vert J_R\cap J_S}$.  Comme $r'_{\vert J_R\cap J_S}=s'_{\vert J_R\cap J_S}$, on en déduit que $r$ et $s$ coïncident sur $J_R\cap J_S$, d'où l'existence de $r+s\in R\bowtie S$, qui vérifie  $(r+s)_{\vert L}=r'+s'$, d'où finalement $r'+s'\in (R\bowtie S)_{\vert L}$.
\begin{flushright}$\square$\end{flushright} 

\subsection{Produit de restrictions}

\begin{prop}\label{prop inclusion produit de restrictions} Pour toute relation multiple $R\in\mathcal{R}$, et pour tout recouvrement $J_R=K\cup L$ du domaine de $R$ par un couple $(L,K)$ de parties de $I$, on a l'inclusion
\[R\subset R_{\vert L}\bowtie R_{\vert K}.\]
\end{prop}
\paragraph{Preuve.} D'après les propositions \ref{prop somme de deux restrictions} et \ref{prop expression du produit de relations par sommes}, on a pour tout $x\in R$ : $x=x_{\vert L}+ x_{\vert K}\in R_{\vert L}\bowtie R_{\vert K}$.
\begin{flushright}$\square$\end{flushright} 

\section{Scissions}

\subsection{Parties propres et bipartitions}

Rappelons qu'une partie propre d'un ensemble $J$ est une partie $K\subset J$ telle que
\[\emptyset\varsubsetneq K\varsubsetneq J.\]

Dans cette section et la suivante, nous ferons souvent appel à des partitions de l'ensemble $I$, ou d'un sous-ensembles $J$ de $I$, constituées de deux parties $K$ et $L$. Nous appellerons de telles partitions des \emph{bipartitions}. Rappelons qu'une partition est un recouvrement constitué de parties propres (donc non \emph{non vides}) et deux à deux disjointes. Ainsi, une bipartition de $J$ est un couple $(K,L)$ de parties de $J$ telles que
\[K\neq \emptyset \neq L\quad\mathrm{et}\quad K\cup L=J \quad\mathrm{et}\quad K\cap L=\emptyset.\] Remarquons que l'existence d'une bipartition de $J$ implique que $J$ a au moins deux éléments.

\subsection{Somme de deux familles de domaines disjoints}

Soient $J$ et $K$ deux parties disjointes  de $I$ : $J\cap K=\emptyset$.
Pour  pour tout couple $(x,y)\in Z_J\times Z_K$, on a $x_{\vert J\cap K}=\bullet=y_{\vert J\cap K}$, donc $x\diamond y$, de sorte que $x+y\in Z_{J\cup K}$ est bien défini.

\subsection{Unicité des factorisations disjointes}

\begin{prop}\label{prop restrictions factorisation} Soient $R\neq 0_{J_R}$ et $S\neq 0_{J_S}$ deux relations non nulles
et de domaines disjoints :
$J_R\cap J_S=\emptyset$.  Alors 
\[R=(R\bowtie S)_{\vert J_R}\quad\mathrm{et}\quad S=(R\bowtie S)_{\vert J_S}.\] 
\end{prop}

\paragraph{Preuve.} D'après la proposition \ref{prop egalite restriction produit}, on a
\[(R\bowtie S)_{\vert J_R}=R\bowtie S_{\vert \emptyset},\]
mais $S\neq 0_{J_S} \Rightarrow S_{\vert \emptyset}=1_\emptyset=\mathbf{1}$,
d'où $(R\bowtie S)_{\vert J_R}=R$. Et de même a-t-on $(R\bowtie S)_{\vert J_S}=S$.
\begin{flushright}$\square$\end{flushright} 

\begin{cor} Deux relations non nulles de domaines disjoints sont nécessairement compatibles.
\end{cor}

\begin{rmq}
La proposition \ref{prop restrictions factorisation} cesse d'être vérifiée si l'on ne suppose pas $J\cap K=\emptyset$ ou si $R$ ou $S$ est nulle. Par exemple, $R\bowtie\mathbf{0} =\mathbf{0}$ et $R\bowtie 0_\emptyset=0_J$ ne déterminent pas $R$. De même, on construit facilement un exemple de relations non nulles incompatibles de domaines respectifs $J$ et $K$ avec $J\cap K\neq \emptyset$, pour lesquelles la proposition n'est évidemment pas vérifiée. 
\end{rmq}

\begin{cor}\label{cor unique factorisation} Soit $T$ une relation non nulle, et soit $(K,L)$ une bipartition de $J_T$. 
Alors, \emph{si elle existe}, une factorisation de $T$ de la forme
\[T=R\bowtie S\] avec $J_R=K$ et $J_S=L$ est nécessairement unique, $R$ et $S$ étant donnés par
\[R=T_{\vert K}\quad\mathrm{et}\quad S=T_{\vert L}.\]
\end{cor}

\subsection{Relations scindables}

\begin{df} Une relation $T$ sera dite \emph{scindable selon une bipartition $(K,L)$ de $J_T$} si elle admet une factorisation de la forme $T=R\bowtie S$ avec $(R,S)\in\mathcal{R}_K\bowtie \mathcal{R}_L$. 
\end{df}

En pratique, les critères suivants  sont extrêmement utiles pour vérifier si une relation est ou non scindable selon une bipartition $(K,L)$.

\begin{prop} \label{prop critere relation scindable} Étant données $T$ une relation dans $\mathcal{E}$ et $(K,L)$ une bipartition de $J_T$, 
 $T$ est scindable selon $(K,L)$  si et seulement si \[T=T_{\vert K}\bowtie T_{\vert L}.\] 
\end{prop}
\paragraph{Preuve.} Si l'égalité a lieu, $T$ est scindable sur $(K,L)$ par définition. Réciproquement, si $T$ est non nulle et scindable selon $(K,L)$, l'égalité résulte immédiatement du corolaire \ref{cor unique factorisation}, tandis que si $T$ est nulle cette égalité est trivialement satisfaite. 
\begin{flushright}$\square$\end{flushright} 

\begin{prop}\label{prop critere somme relation scindable} Une relation multiple $T$ est scindable selon une bipartition $(K,L)$ de $J_T$ si et seulement si on a 
 \[\forall x\in T_{\vert K}, \forall y\in T_{\vert L}, x+y\in T.\]
\end{prop}
\paragraph{Preuve.} D'après les propositions \ref{prop inclusion produit de restrictions} et \ref{prop critere relation scindable}, $T$ est scindable selon $(K,L)$ si et seulement si on a l'inclusion $T_{\vert K}\bowtie T_{\vert L}\subset T$, et la proposition \ref{prop expression du produit de relations par sommes} achève la preuve.
\begin{flushright}$\square$\end{flushright}

\begin{df} Une relation $T$ est dite \emph{scindable} s'il existe une bipartition $(K,L)$ de $J_T$ telle que $T$ soit scindable selon $(K,L)$.
\end{df}

\begin{rmq} Puisque l'existence d'une bipartition de $J_T$ implique que $J_T$ a au moins deux éléments, aucune relation de domaine vide ou singleton n'est scindable.
\end{rmq}

\subsection{Parties scindables pour une relation $R$}

\begin{df} Étant donnée $R\in\mathcal{R}$ 
une relation multiple sur $\mathcal{E}$, 
une partie $J$ de $J_R$ est dite \emph{scindable pour $R$} 
si la relation $R_{\vert J}$ est scindable.
\end{df}

D'après la proposition \ref{prop critere relation scindable}, $J\subset J_R$ est scindable pour $R$ s'il existe deux parties non vides complémentaires $K$ et $L$ dans $J$ --- ce qui suppose que $J$ soit de cardinal $\geq 2$ --- telles que
\begin{equation}\label{equ critere partie scindable}
R_{\vert J}=R_{\vert K}\bowtie R_{\vert L}.
\end{equation}

\subsection{Parties détachables d'une relation $R$}

Dans cette section, une relation  $R\in \mathcal{R}_\mathcal{E}$ est donnée et, pour tout $J\subset J_R$, on pose  \[\neg J=J_R\setminus J.\]

\begin{df} Une partie $J\subset J_R$ est dite \emph{détachable de $R$} si l'on a
\[R=R_{\vert \neg J}\bowtie 1_J.\]
\end{df}

Dans le cas où $J$ est une partie propre non vide de $J_R$, autrement dit lorsque $\emptyset\varsubsetneq J\varsubsetneq J_R$, le couple $(J,\neg J)$ est une bipartition de $J_R$, et la définition précédente revient à dire que $R$ est scindable selon cette bipartition, avec en outre \[R_{\vert J}=1_J.\] Bien entendu, cette dernière condition n'est pas suffisante pour faire de $J$ une partie détachable de $R$.

\begin{exm}
Pour $I=\{1,2\}$, $J=\{1\}$ et $f$ une application quelconque $E_1\to E_2$ de graphe $G$, la relation $R=(I,G)$ vérifie nécessairement $R_{\vert J}=1_J$. En outre, $J$ est détachable de $R$ si et seulement si $f$ est constante.
\end{exm}

\begin{lm}\label{lm detachable} $J$ est une partie de $J_R$ détachable de $R$ si et seulement si on a \emph{pour tout} $x\in Z_{J_R}$
\[\rho_{(J_R,\neg J)} (x) \in \rho_{(J_R,\neg J)} (R)\Rightarrow x\in R.\]
\end{lm}

\paragraph{Preuve.} Par définition du produit $\rho_{(J_R,\neg J)}(R)\bowtie 1_J$, on a $J$ détachable de $R$ si et seulement si
\[
G_R= \rho_{(J_R,\neg J)}^{-1}(\rho_{(J_R,\neg J)}(G_R))\cap \rho_{(J_R,J)}^{-1}(Z_J)
\linebreak
=\rho_{(J_R,\neg J)}^{-1}(\rho_{(J_R,\neg J)}(G_R)).
\]
Puisque l'inclusion $G_R\subset \rho_{(J_R,\neg J)}^{-1}(\rho_{(J_R,\neg J)}(G_R))$ est toujours trivialement satisfaite, $J$ est donc détachable de $R$ si et seulement si on a l'inclusion réciproque
\[\rho_{(J_R,\neg J)}^{-1}(\rho_{(J_R,\neg J)}(G_R))\subset G_R,\]
autrement dit si pour tout $x\in Z_{J_R}$ tel que $\rho_{(J_R,\neg J)}(x)\in \rho_{(J_R,\neg J)}(G_R)$, on a $x\in G_R$.
\begin{flushright}$\square$\end{flushright} 

\begin{prop} Si $J$ est détachable de $R$, alors pour toute $J_R$-famille $R$-compatible $y$ et pour toute $J$-famille $x$, on a
\[\rho_{(J_R,\neg J)}(y)+x\in R.\]
\end{prop}
\paragraph{Preuve.} Posons $z=\rho_{(J_R,\neg J)}(y)+x$. On a $\rho_{(J_R,\neg J)}(z)=\rho_{(J_R,\neg J)}(y)\in \rho_{(J_R,\neg J)}(R)$, donc $z\in R$.
\begin{flushright}$\square$\end{flushright} 

\begin{prop} Si $J$ et $K$ sont deux parties de $J_R$ détachables de la relation $R\in \mathcal{R}_\mathcal{E}$, alors $J\cup K$ est également une partie détachable de $R$.
\end{prop}

\paragraph{Preuve.}  Posons $C=\neg(J\cup K)=\neg J\cap\neg K$. Soit $x\in Z_{J_R}$ quelconque tel que $\rho_{(J_R,C)} (x) \in \rho_{(J_R,C)} (R)$. 
D'après le lemme \ref{lm detachable}, il suffit de prouver  que $x\in R$. Posons $x_C=\rho_{(J_R,C)} (x)$. Puisque par hypothèse $x_C\in \rho_{(J_R,C)}(R)$, il existe $y\in R$ tel que $x_C=\rho_{(J_R,C)}(y)$. 
Puisque $J$ est détachable, on déduit de $\rho_{(J_R,\neg J)}(y)\in \rho_{(J_R,\neg J)}(R)$ que $z=\rho_{(J_R,\neg J)}(y)+\rho_{(J_R,J)}(x)$ est également $R$-compatible. 
Par conséquent $\rho_{(J_R,\neg K)}(z)\in \rho_{(J_R,\neg K)}(R)$, et $K$ étant détachable on en déduit que $w=\rho_{(J_R,\neg K)}(z)+\rho_{(J_R, K)}(x)$ est lui aussi $R$-compatible. Or, par construction même, $w$ et $x$ coïncident sur $K$, sur $J\cap \neg K$ et sur $C$, de sorte que $w=x$, d'où $x\in R$.
\begin{flushright}$\square$\end{flushright} 

\begin{rmq}\label{rmq contre exemple fluide} L'union d'une famille finie de parties de $J_R$ détachables de $R$ est donc encore détachable de $R$. Par contre, ceci n'est pas vrai en général pour une famille quelconque. Considérons par exemple le cas où $I=\mathbf{N}$, $E_i=\mathbf{N}$ pour tout $i\in I$, et $R=(I,G)$ avec $G$ l'ensemble des suites $x\in \mathbf{N}^\mathbf{N}$ comportant une infinité de zéros : $x\in G \Leftrightarrow Card(\{n\in\mathbf{N},x_n=0\})=\aleph_0$. Alors $\{n\}$ est détachable de $R$ pour tout $n\in\mathbf{N}$, mais $I=\mathbf{N}$ lui-même n'est pas détachable de $R$ puisque $R$ n'est pas la relation triviale.
\end{rmq}

\subsubsection{Partie externe et socle d'une relation}

\begin{df} On appelle \emph{partie externe $Ex(R)$ d'une relation $R$} l'union des parties de $I$ qui sont détachables de $R$. On appelle \emph{socle de $R$} l'ensemble $Soc(R)=J_R\setminus Ex(R)$. 
\end{df}

\begin{exm}
Toute relation triviale $1_J$ a un socle vide.
\end{exm}

\begin{df}
Une relation sera dite 
\begin{itemize}
\item \emph{mouvante} si sa partie externe est non détachable,
\item \emph{ancrée} si elle n'est pas mouvante,
\item \emph{fluide} si elle n'est pas triviale mais que son socle est vide,
\item \emph{solide} si sa partie externe est vide.
\end{itemize}
\end{df}

Une relation fluide est nécessairement mouvante, tandis qu'une relation solide est nécessairement ancrée. Une relation ancrée est déterminée par sa restriction à son socle, tandis qu'il ne suffit pas de connaître une relation mouvante sur son socle pour la connaître entièrement. Les relations finies, sans être nécessairement solides, sont toujours ancrées. 

\begin{exm} L'exemple de la remarque \ref{rmq contre exemple fluide} est celui d'une relation fluide.
\end{exm}

\section{Structure connective d'une relation $R\in\mathcal{R}$}

Rappelons\footnote{Voir \cite{Dugowson:201012} et \cite{Dugowson:201203}.  } qu'un \emph{espace connectif} $(X,\mathcal{K})$ est la donnée d'un ensemble de points $X$, appelé support de l'espace, et d'un ensemble $\mathcal{K}$ de parties de $X$, appelé \emph{structure connective} de l'espace, tel que
\begin{displaymath}
\forall \mathcal{I}\in \mathcal{P}(\mathcal{K}), \left(\bigcap_{K\in\mathcal{I}}K\ne\emptyset\Rightarrow \bigcup_{K\in\mathcal{I}}K\in\mathcal{K}\right).
\end{displaymath} 
Notons que la propriété ci-dessus entraîne en particulier que  $\emptyset\in\mathcal{K}$. L'espace connectif $(X,\mathcal{K})$ est dit \emph{intègre} si les singletons $\{x\}$, où $x\in X$, appartiennent tous à $\mathcal{K}$. Pour tout ensemble de parties $\mathcal{C}\subset \mathcal{P}X$, on note 
\[[\mathcal{C}]\]
la structure connective engendrée par $\mathcal{C}$.

Pour démontrer le théorème \ref{thm Brunn}, nous aurons besoin du lemme suivant.

\begin{lm}\label{lm partie non connexe} Étant donné $(X,\mathcal{K})$ un espace connectif intègre, et $A\subset X
$ une partie \emph{non connexe} de $X$, il existe nécessairement une bipartition $(L,M)$ de $A$ telle que pour toute partie connexe $K\in\mathcal{K}$ incluse dans $A$ on a
\[
\left\vert
\begin{array}{l}
\mathrm{ou}\,\,\mathrm{bien}\quad K\subset L,\\
\mathrm{ou}\,\,\mathrm{bien}\quad K\subset M.\\
\end{array}
\right.
\]
\end{lm}
\paragraph{Preuve.} $A$ étant supposé non connexe est nécessairement non vide. Soit $x\in A$. Prenons pour $L$ la composante connexe\footnote{Voir \cite{Dugowson:201012}.} de $x$ dans l'espace connectif induit par $(X,\mathcal{K})$ sur $A$, c'est-à-dire dans l'espace $(A,\mathcal{K}\cap\mathcal{P}A)$. Autrement dit, $L$ est le plus grand connexe inclus dans $A$ et contenant $x$. Puisque $A$ est non connexe, on a nécessairement $\emptyset\varsubsetneq L\varsubsetneq A$, de sorte que $(L,M)$ forme une bipartition de $A$, où $M=A\setminus L$. 
Soit maintenant $K$ une partie connexe incluse dans $A$. On a soit $K\cap L=\emptyset$, et dans ce cas $K\subset M$, soit $K\cap L\neq\emptyset$ et dans ce cas $K\cup L$ est un connexe contenant $x$ et contenu dans $A$, de sorte que $K\cup L\subset L$, autrement dit $K\subset L$.
\begin{flushright}$\square$\end{flushright}

\subsection{Définition}

\begin{prop} Étant donnée $R$ une relation dans $\mathcal{E}$, l'ensemble $\mathcal{K}_R$ des parties de $J_R$ non scindables pour $R$ constitue une structure connective intègre sur $J_R$.
\end{prop}
\paragraph{Preuve.} La partie vide et les singletons font nécessairement partie de $\mathcal{K}_R$, puisqu'une partie scindable pour $R$, admettant une bipartition, a nécessairement au moins deux éléments. Soit maintenant $\mathcal{C}\subset \mathcal{K}_R$ un ensemble de parties non scindables pour $R$ tel que 
\[\bigcap_{C\in\mathcal{C}}C\neq\emptyset.\]  Montrons par l'absurde que $U=\bigcup_{C\in\mathcal{C}}C$ est également non scindable. Si $U$ était scindable pour $R$, il admettrait une bipartition $(L,M)$ telle que $R_{\vert U}$ soit scindable selon $(L,M)$, de sorte que l'on aurait d'après la proposition \ref{prop critere relation scindable}
\[R_{\vert U}=R_{\vert L}\bowtie R_{\vert M}.\]
Soit $x\in\displaystyle\bigcap_{C\in\mathcal{C}}C$. On a soit $x\in L$, soit $x\in M$. Supposons pour fixer les idées que $x\in L$. Puisque $M\neq \emptyset$, il existe $C\in \mathcal{C}$ tel que $M\cap C\neq \emptyset$. Mais on a aussi $L\cap C\neq \emptyset$, puisque $x\in L\cap C$. Donc $(L\cap C, M\cap C)$ est une bipartition de $C$. Mais d'après la proposition \ref{prop egalite restriction produit}, on aurait 
\[R_{\vert C}=R_{\vert L\cap C}\bowtie R_{\vert M\cap C},\] de sorte que $C$ serait scindable pour $R$, ce qui est absurde.
\begin{flushright}$\square$\end{flushright} 

\begin{df} On appelle structure connective d'une relation multiple $R\in\mathcal{R}$ la structure connective $\mathcal{K}_R$ définie dans la proposition précédente.
\end{df}

\subsection{Exemple}

On pourrait penser que, pour tout couple de relations multiples $(R,S)\in\mathcal{R}^2$, la structure connective de $R\bowtie S$ devrait être incluse dans la structure connective engendrée par les connexes de $R$ et ceux de $S$. En général, cela est faux, de sorte que
\[\mathcal{K}_{R\bowtie S}\nsubset [\mathcal{K}_R\cup \mathcal{K}_S].\]
Donnons-en un contre exemple simple dans le cas où $I=\{1,2,3\}$ et, pour tout $i\in I$, $E_i=\{0,1\}$. On définit la relation $R$ de la façon suivante
\[(x_1,x_2,x_3)\in R \Leftrightarrow  \exists i\in I, x_i=0,\] et la relation $S$ par
\[(x_1,x_2,x_3)\in S \Leftrightarrow  \exists i\in I, x_i=1.\] On vérifie facilement que $\mathcal{K}_R=\mathcal{K}_S=\mathcal{B}_3$, la structure borroméenne\footnote{Voir \cite{Dugowson:201012}.} sur $I$, tandis que $\mathcal{K}_{R\bowtie S}$ est la structure connective grossière que $I$.


\subsection{Théorème de Brunn}

%
%

En référence au résultat annoncé par Brunn\footnote{Mais non entièrement démontré par lui, puisqu'il faudra attendre Kanenobu\cite{Kanenobu:198504} en 1984 pour avoir une telle démonstration complète.} \cite{Brunn:1892a} en 1892 à propos de la structure des entrelacs, j'appelle \og théorème de Brunn\fg\, relatif à une classe d'objets à chacun desquels se trouve associée une structure connective l'énoncé affirmant que toute structure connective --- ou du moins toute structure connective d'un certain type, par exemple toute structure connective intègre finie --- est celle d'au moins un objet de cette classe.

\begin{thm}\label{thm Brunn}  Pour tout ensemble $I$, il existe un choix des ensembles $E_i$ tel que pour toute structure connective intègre $\mathcal{K}$ sur $I$ il existe une relation $R\in\mathcal{R}$ telle que $\mathcal{K}_R=\mathcal{K}$.
\end{thm}

\paragraph{Preuve.} Considérons la construction suivante. On prend le même ensemble $E_i$ pour tous les $i\in I$, à savoir $E_i=\{0,1\}^{\mathcal{P}(I)}$, ensemble des applications de l'ensemble des parties de $I$ dans $\{0,1\}$. Soit maintenant  $\mathcal{K}$ une structure connective intègre sur $I$. Considérons la relation multiple $R$ dans $\mathcal{E}$, de domaine $I$, définie de la façon suivante : une famille $(f_i:\mathcal{P}(I)\longrightarrow \{0,1\})_{i\in I}$ est  $R$-compatible si et seulement si
\[
\forall K\in\mathcal{K}\setminus \{\emptyset\}, \exists i\in K, f_i(K)=1.
\]
et vérifions que la structure connective de $R$ est précisément $\mathcal{K}$. 

Avant toute chose, commençons par remarquer que pour toute partie $A\subset I$, la restriction $R_{\vert A}$ de $R$ à $A$ a pour graphe l'ensemble des $A$-familles $(f_i)_{i\in A}$ telles que
\begin{equation}\label{eqn caracterisation de R restreinte a A}
\forall K\in\mathcal{K}\cap \mathcal{P}A\setminus \{\emptyset\}, \exists i\in K, f_i(K)=1.
\end{equation} En effet, cette condition devant être satisfaite par toute $I$-famille $R$-compatible doit également, par restriction, être satisfaite par toute $A$-famille $R_{\vert A}$-compatible. Réciproquement, si une $A$-famille vérifie la condition en question, il est aisé de la prolonger en une $I$-famille $R$-compatible, puisqu'il suffit pour tout $i\in I\setminus A$ et toute partie $B\subset I$ de poser $f_i(B)=1$ pour obtenir un tel prolongement.

Si $K\in\mathcal{K}$, alors $R$ ne peut pas être scindable sur $K$. Raisonnons par l'absurde : supposons que $R_{\vert K}$ soit scindable selon une bipartition $(L,M)$ de $K$. Alors la famille \linebreak $l=(l_i:\mathcal{P}(I)\longrightarrow \{0,1\})_{i\in L}$ définie pour tout $i\in L$ et toute partie $A\subset I$ par 
\[
\left\{
\begin{array}{l}
l_i(A)=1\quad\mathrm{si}\,A\neq K,\\
l_i(K)=0
\end{array}
\right.
\] 
est $R_{\vert L}$-compatible,
 car en la prolongeant sur $I$ par la famille $\tilde{l}$ définie pour $i\in L$ par $\tilde{l}_i=l_i$ et pour $i\in I\setminus L$ par $\tilde{l}_i(A)=1$ pour \emph{toute} partie $A\subset I$ (y compris, donc, pour $A=K$), on obtient $\tilde{l}\in R$ puisque la propriété caractérisant la relation $R$ est trivialement satisfaite.
De même, la famille  $m=$ $(m_i:\mathcal{P}(I)\longrightarrow \{0,1\})_{i\in M}$ définie pour tout $i\in M$ et toute partie $A\subset I$ par 
\[
\left\{
\begin{array}{l}
m_i(A)=1\quad\mathrm{si}\,A\neq K,\\
m_i(K)=0
\end{array}
\right.
\] est $R_{\vert M}$-compatible. Par conséquent, $l+m$ doit être $R_{\vert K}$ compatible, ce qui est absurde car il n'existe pas d'indice $i\in K$ tel que $(l+m)_i(K)=1$.

Réciproquement, soit $K\subset I$ non scindable pour $R$. Montrons que $K\in\mathcal{K}$. Nous allons à nouveau raisonner par l'absurde. Supposons que $K\notin\mathcal{K}$. Dans cette hypothèse, d'après le lemme \ref{lm partie non connexe}, il doit exister une bipartition $(L,M)$ de $K$ telle que toute partie connexe $C\in\mathcal{K}$ incluse dans $K$ vérifie soit $C\subset L$, soit $C\subset M$. Soit alors $f\in R_{\vert L}$ et $g\in R_{\vert M}$. La somme $f+g$ est une $K$-famille qui est nécessairement $R_{\vert K}$-compatible, puisque pour tout connexe $C\subset K$, on a soit $C\subset L$, d'où l'existence de $i\in C$ tel que $1=f_i(C)=(f+g)_i(C)$, soit $C\subset M$, auquel cas il existe $i\in C$ tel que $1=g_i(C)=(f+g)_i(C)$, de sorte que la relation (\ref{eqn caracterisation de R restreinte a A}) est satisfaite. Par conséquent, d'après la proposition \ref{prop critere somme relation scindable}, la relation $R$ est scindable sur $(L,M)$, ce qui contredit l'hypothèse qui avait été faite. 

Finalement, nous avons établi que la structure connective de la relation $R$ ainsi définie est bien la structure donnée $\mathcal{K}$, ce qui prouve le théorème.
\begin{flushright}$\square$\end{flushright} 

\bibliographystyle{plain}




\newpage

\tableofcontents

\end{document}